\documentclass[12pt]{article}
\usepackage{fullpage,graphicx,psfrag,amsmath,amsfonts,verbatim}
\usepackage[small,bf]{caption}
\usepackage{algorithmicx,algpseudocode}
\usepackage{url,algorithm,amsmath,amsthm}
\usepackage{siunitx}
\newcommand{\BEAS}{\begin{eqnarray*}}
\newcommand{\EEAS}{\end{eqnarray*}}
\newcommand{\BEA}{\begin{eqnarray}}
\newcommand{\EEA}{\end{eqnarray}}
\newcommand{\BEQ}{\begin{equation}}
\newcommand{\EEQ}{\end{equation}}
\newcommand{\BIT}{\begin{itemize}}
\newcommand{\EIT}{\end{itemize}}
\newcommand{\BNUM}{\begin{enumerate}}
\newcommand{\ENUM}{\end{enumerate}}

\newcommand{\BA}{\begin{array}}
\newcommand{\EA}{\end{array}}

\newcommand{\eg}{{\it e.g.}}
\newcommand{\ie}{{\it i.e.}}

\newcommand{\reals}{{\mbox{\bf R}}}

\newcommand{\Co}{\mathop {\bf Co}}

\newcommand{\argmin}{\mathop{\rm argmin}}

{\begin{quote}}{\end{quote}}

\makeatletter
\long\def\@makecaption#1#2{
   \vskip 9pt 
   \begin{small}
   \setbox\@tempboxa\hbox{{\bf #1:} #2}
   \ifdim \wd\@tempboxa > 5.5in
        \begin{center}
        \begin{minipage}[t]{5.5in}
        \addtolength{\baselineskip}{-0.95pt}
        {\bf #1:} #2 \par
        \addtolength{\baselineskip}{0.95pt}
        \end{minipage}
        \end{center}
   \else 
	\hbox to\hsize{\hfil\box\@tempboxa\hfil}  
   \fi
   \end{small}\par
}
\makeatother

\newcounter{oursection}

\newcounter{lecture}

\usepackage{subcaption}
\usepackage{circuitikz}
\usepackage{tikz}

\title{A Simple Effective Heuristic for Embedded Mixed-Integer 
Quadratic Programming}
\author{Reza Takapoui \and Nicholas Moehle 
\and Stephen Boyd \and Alberto Bemporad}

\begin{document}
\maketitle

\begin{abstract}
In this paper we propose a fast optimization algorithm 
for approximately minimizing convex quadratic functions
over the intersection of affine and separable constraints
(\ie, the Cartesian product of possibly nonconvex real sets).
This problem class contains many NP-hard problems
such as mixed-integer quadratic programming.
Our heuristic is based on a variation of the 
alternating direction method of multipliers (ADMM),
an algorithm for solving convex optimization problems.
We discuss the favorable computational aspects of our algorithm,
which allow it to run quickly even on very modest computational
platforms such as embedded processors.
We give several examples for which an approximate
solution should be found very quickly,
such as management of a hybrid-electric vehicle drivetrain
and control of switched-mode power converters.
Our numerical experiments suggest that our method
is very effective in finding a feasible point with small objective value;
indeed, we find that in many cases, it finds the global solution.
\end{abstract}

\section{Introduction}

\subsection{The problem}
We consider the problem
\BEQ
\label{miqp}
\begin{array}{ll}
\mbox{minimize}   & (1/2)x^T P x + q^T x + r\\
\mbox{subject to} & Ax = b\\
& x \in \mathcal X
\end{array}
\EEQ
with decision variable $x\in \reals^n$.
The problem parameters are 
the symmetric positive semidefinite matrix $P\in\reals^{n\times n}$,
the matrix $A \in \reals^{m\times n}$,
the vectors $b \in \reals^m$ and $q \in \reals^n$, and the real number $r\in \reals$.
The constraint set $\mathcal X$ is the Cartesian product of
(possibly nonconvex) real, closed, nonempty sets, \ie, $\mathcal X = 
\mathcal X_1 \times \cdots \times \mathcal X_n$, where $\mathcal X_i
\subseteq \reals$ are closed, nonempty 
subsets of $\reals$ for $i=1,\ldots,n$.
If $\mathcal X_i$ is a convex set, we refer to variable $x_i$ as a
\emph{convex variable},
and if $\mathcal X_i$ is a nonconvex set, we call variable 
$x_i$ a \emph{nonconvex variable}.

\paragraph{Some applications.}
Many problems can be put into the form of problem (\ref{miqp}).
For example, if some of the sets $\mathcal X_i$ are subsets
of integers, our formulation addresses 
mixed-integer quadratic and mixed-integer linear programs.
This includes applications such as 
admission control \cite{oulai2007new}, 
economic dispatch \cite{papageorgiou2007mixed}, 
scheduling \cite{catalao2010scheduling}, 
hybrid vehicle control \cite{murgovski2012component},
thermal unit commitment problems \cite{carrion2006computationally}, 
Boolean satisfiability problems (SAT) \cite{jeroslow1990solving},
and hybrid model predictive control \cite{bemporad1999control}.
Another application is embedded signal decoding in communication
systems, when the nonconvex sets are signal constellations 
(\eg, QAM constellations; see \cite[pg. 416]{glover2010digital}).

\paragraph{Complexity.}
If $\mathcal X$ is a convex set, problem \eqref{miqp} is a convex
optimization problem and can be readily solved using standard
convex optimization techniques. Otherwise,
the problem \eqref{miqp} can be hard in general. 
It trivially generalizes mixed-integer quadratic programming,
an NP-complete problem, and can therefore be used to encode
other NP-complete problems such as the traveling salesman problem
(TSP) \cite{papadimitriou1998combinatorial}, 
Boolean satisfiability (SAT) \cite{li2004satisfiability,karp1972reducibility},
set cover \cite{hochbaum1982approximation}, 
and set packing \cite{padberg1973facial}. 
Hence, any algorithm 
that guarantees finding the global solution to \eqref{miqp} 
suffers from non-polynomial worst-case time
(unless $\mathrm {P = NP}$).

\subsection{Solve techniques}
\paragraph{Exact methods.}
There are a variety of methods for solving \eqref{miqp} exactly.
When all of the nonconvex sets $\mathcal X_i$ in \eqref{miqp} are
finite, the simplest method is brute force; enumerating through all
possible combinations of discrete variables and solve a convex
optimization problem for each possible combination and finding
the point with the smallest objective value. Other methods
such as branch-and-bound \cite{lawler1966branch} and 
branch-and-cut \cite{stubbs1999branch}
are guaranteed to find the global solution.
Cutting plane methods \cite{gomory1958outline, chvatal1989cutting} 
rely on solving the relaxation and adding a linear constraint to
drive the solution towards being integer.
Special purpose methods have been introduced for some 
specific subclasses of \eqref{miqp}.
Unfortunately, these methods have non-polynomial 
worst-case runtime, and are often burdensome to use in practice,
especially for embedded optimization,
where runtime, memory limits, and code simplicity are prioritized.
Also, these methods suffer from a large variance in the
algorithm runtime.

\paragraph{Heuristics.}
On the other hand, many heuristics have been introduced that
can deliver a good, but suboptimal (and possibly infeasible)
point in a very short amount of time.
For example, the \emph{relax-and-round} heuristic
consists of replacing each $\mathcal X_i$ by its convex hull,
solving the resulting relaxation (a convex quadratic program),
and projecting the solution onto the nonconvex constraint sets.
Another heuristic is to fix the nonconvex variables for several 
reasonable guesses and solve the convex optimization problem 
for convex variables. 
(Each of these problems may not find a feasible point, even if one exists.)
The \emph{feasibility pump} is a heuristic
to find a feasible solution to a generic mixed integer program 
and is discussed in \cite{fischetti2005feasibility, 
bertacco2007feasibility, achterberg2007improving}.
Such heuristics are often quite effective,
and can be implemented on very modest computational hardware,
making them very attractive for embedded applications
(even without any theoretical guarantees).

\subsection{Embedded applications}
We focus on embedded applications where 
finding a feasible point with relatively small objective will often result
in performance that is practically indistinguishable from 
implementing the global solution. In embedded applications,
the computational resources are limited and a solution must
be found in a small time. Hence, methods to find the global solution
are not favorable, because their large variance in runtime cannot be
tolerated.

In an embedded application,
it is often required to solve several instances of (\ref{miqp}),
with different values of the parameters.
Here we distinguish two separate use cases,
depending on whether one or both of $P$ or $A$ change.
This distinction will play an important role in solution methods.
In the first use case, we solve many instances of (\ref{miqp})
in which any of the parameters may change between instances.
In the second use case, we solve instances of (\ref{miqp})
in which $q$, $b$, and $\mathcal X$ change between instances,
but $P$ and $A$ are constant.
Although this is more restrictive than the first use case,
many applications can be well modeled using this approach,
including linear, time-invariant model predictive control
and moving horizon estimation.
Indeed, all of the three examples we present in \S\ref{s-examples}
are of this type.

\subsection{Contributions}
Our proposed algorithm is a simple and computationally efficient
heuristic to find approximate
solutions to problem \eqref{miqp} quickly.
It is based on the alternating direction method of multipliers 
(ADMM), an algorithm for solving convex optimization problems.
Because the problem class we address 
includes nonconvex optimization problems, our method is not 
guaranteed to find the global solution, or even converge.

Numerical experiments suggest that this heuristic is an effective
tool to find the global solution in a variety of problem instances. 
Even if our method does not find the global solution,
it usually finds a feasible point with reasonable objective value.
This makes it effective for many
embedded optimization applications, where
finding a feasible point with relatively small objective value often results
in performance that is practically indistinguishable from 
implementing the global solution.

Comparison of the runtime with commercial solvers such as 
MOSEK \cite{mosek} and CPLEX \cite{cplex2009v12} show 
that our method can be
substantially faster than solving a global optimization method,
while having a competitive practical performance.

\subsection{Related work}

\paragraph{Fast embedded optimization.}
In recent years,
much research has been devoted to solving moderately-sized 
convex optimization problems
quickly (\ie, in milliseconds or microseconds), 
possibly on embedded platforms.
Examples include the SOCP solvers 
ECOS \cite{domahidi2013ecos},
and FiordOs \cite{ullmann2011fiordos},
and the QP solver
CVXGEN \cite{mattingley2012cvxgen}.
Other algorithms have been developed exclusively for convex optimal 
control problems; see 
\cite{wang2010fast,odonoghue2012splitting,jerez2014embedded}.
In addition, recent advances in automatic code generation for convex 
optimization \cite{mattingley2011receding,chu2013code}
can significantly reduce the cost
and complexity of using an embedded solver. 
Some recent effort has been devoted to (globally) 
solving mixed-integer convex programs
very quickly;
see \cite{bemporad2015solving}, 
\cite{frick2015embedded} and references therein.

\paragraph{Nonconvex ADMM.}
Even though ADMM was originally introduced as a tool for convex
optimization problems, it turns out to be a powerful heuristic method 
even for NP-hard nonconvex problems 
\cite[\S 5, 9]{boyd2011distributed}. Recently, this tool has been 
used as a heuristic to find approximate solutions to nonconvex 
problems \cite{chartrand2013nonconvex,chartrand2012nonconvex}.
In \cite{derbinsky2013improved}, the authors study the
\emph{Divide and Concur} algorithm as a special case of a message-passing version of the 
ADMM, and introduce a three weight version of this algorithm
which greatly improves the performance for some nonconvex 
problems such as circle packing and the Sudoku puzzle.

\section{Our heuristic}
\subsection{Algorithm}
Our proposed algorithm is an extension of the alternating direction 
method of multipliers (ADMM) for constrained optimization to the
nonconvex setting \cite[\S 5,9]{boyd2011distributed}. 
ADMM was originally introduced for solving convex
problems, but practical evidence suggests that it can be an effective 
method to approximately solve some nonconvex problems as well.
In order to use ADMM, we rewrite problem \eqref{miqp} as 
\BEQ\label{miqp_consensus}
\begin{array}{ll}
\mbox{minimize}   & (1/2)x^T P x + q^T x + I_{\mathcal X}(z)\\
\mbox{subject to} &  \left[ \begin{array}{c}
A \\
I \end{array} \right]x - 
\left[ \begin{array}{c}
0 \\
I \end{array} \right]z = 
\left[ \begin{array}{c}
b \\
0 \end{array} \right]. \\
\end{array}
\EEQ
Here $I_{\mathcal X}$ denotes the indicator function of $\mathcal X$, 
so that $I_{\mathcal X}(x)=0$ for $x\in\mathcal X$ and 
$I_{\mathcal X}(x)=\infty$ for $x\notin\mathcal X$. 
Each iteration in the algorithm consists of the 
following three steps:
\BEAS
x^{k+1/2} &:= & \argmin_x \left((1/2)x^T P x + q^T x+ (\rho/2)\left\|\left[ \begin{array}{c}
A \\
I \end{array} \right]x - 
\left[ \begin{array}{c}
0 \\
I \end{array} \right]
x^k -
\left[ \begin{array}{c}
b \\
0 \end{array} \right]
+ u^k \right\|_2^2\right)\\
x^{k+1} &:= & \Pi\left(x^{k+1/2} + 
\left[ \begin{array}{cc}
0 & I \end{array} \right]
u^k\right)\\
u^{k+1} &:= & u^k + \left[ \begin{array}{c}
A \\
I \end{array} \right]x^{k+1/2} - 
\left[ \begin{array}{c}
0 \\
I \end{array} \right]x^k - 
\left[ \begin{array}{c}
b \\
0 \end{array} \right].
\EEAS
Here, $\Pi$ denotes the projection onto $\mathcal X$. 
Note that if $\mathcal X$ is not convex, 
the projection onto $\mathcal X$ may not be unique;
for our purposes, we only need that
$\Pi(z) \in\argmin_{x\in \mathcal X} \|x - z\|_2$
for all $z\in \reals^n$. Since $\mathcal X$ is the Cartesian 
product of subsets of the real line, \ie,
$\mathcal X = \mathcal X_1 \times \cdots \times \mathcal X_n$,
we can take $\Pi(z) = \Pi_1(z_1) \times \cdots \times \Pi_n(z_n)$,
where $\Pi_i$ is a projection function onto $\mathcal X_i$.
Usually evaluating $\Pi_i(z)$ is inexpensive; 
for example, if $\mathcal X_i=[a,b]$ is an interval,
$\Pi_i(z) = \min\{\max\{z,a\},b\}$. If $\mathcal X_i$ is the set of 
integers, $\Pi_i$ rounds its argument to the nearest integer. 
For any finite set $\mathcal X_i$ with $k$ elements,
$\Pi_i(z)$ is a closest point to $z$ that belongs to 
$\mathcal X_i$, which can be found by $\lceil \log_2k \rceil$ comparisons.

\subsection{Convergence} 
If the set $\mathcal X$ is convex and problem \eqref{miqp} is feasible, 
the algorithm is guaranteed to converge to an optimal point
\cite[\S3]{boyd2011distributed}. 
However, for $\mathcal X$ nonconvex, there is no such guarantee. 
Indeed, because problem \eqref{miqp} can be NP-hard, 
any algorithm that finds the global solution 
suffers from nonpolynomial worst-case runtime. 
Our approach is to give up the accuracy and use 
methods that find an approximate solution in a small time. 

Our numerical results verify that even for simple examples,
the algorithm may fail to converge, converge to a suboptimal point, 
or fail to find a feasible point, even if one exists. 
Since 
the objective value need not decrease monotonically (or at all),
it is critical to keep track of the best point found runtime. 
That is,
for a selected primal feasibility tolerance $\epsilon^{\mathrm{tol}}$, 
we shall reject all points $x$ such that 
$\|Ax - b\| > \epsilon^{\mathrm{tol}}$, and among those 
primal feasible points $x$ that $\|Ax - b\| \leq \epsilon^{\mathrm{tol}}$, 
we choose the point with the smallest objective
value. Here, $\epsilon^{\mathrm{tol}}$ is a tolerance for accepted feasibility.
We should remind the reader again, that this point need 
not be the global minimum.

\subsection {Initialization} 
To initialize $x^{0}$, one can randomly choose a point in 
$\Co \mathcal X$,
where $\Co \mathcal X$ denotes the convex hull of $\mathcal X$.
More specifically, this means that we need to have access to a 
subroutine that generates random points in $\Co \mathcal X$.
Our numerical results show that running the algorithm multiple 
times with different random initializations increases the chance of 
finding a feasible point with smaller objective value. Hence, we
suggest running the algorithm multiple times initialized with random
starting points and report the best point as the approximate solution.
We always initialize $u^{0}=0$.

\subsection{Computational cost}
In this subsection, we make a few comments about the computational
cost of each iteration.
The first step involves minimizing a strongly convex quadratic function
and is actually a linear operator. 
The point $x^{k+1/2}$ can be found by solving the following 
system of equations:
\[ 
\left[ \begin{array}{cc}
P+\rho I & A^T \\
A & -(1/\rho) I
\end{array} \right]
\left[ \begin{array}{cc}
x^{k+1/2} \\
v 
\end{array} \right]=
\left[ \begin{array}{cc}
-q+ \rho \left( x^k + A^Tb - \left[ \begin{array}{cc}
A^T & I\end{array} \right]u^k \right)\\
0
\end{array} \right].
\] 
Since the matrix on the lefthand side remains constant for all iterations,
we can precompute the $LDL^T$ factorization of this matrix 
once and cache the factorization for use in 
subsequent iterations. When $P$ and $A$ are dense, 
the factorization cost is $O(n^3)$
yet each subsequent iteration costs only $O(n^2)$.
(Both factorization and solve costs 
can significantly smaller if $P$ or $A$ is sparse.)
Amortizing the factorization step over all iterations
means that the first step is quite efficient.
Also notice that the matrix on the lefthand side is 
quasi-definite and hence favorable for $LDL^T$ factorization.

In many applications, $P$ and $A$ do not change across problem 
instances.
In this case, for different problem instances, we solve 
\eqref{miqp} for the same $P$ and $A$ and varying $b$ and $q$. 
This lets us use the same $LDL^T$ factorization,
which results in a significant saving in computation.

The second step involves projection onto 
$\mathcal X = \mathcal X_1 \times \cdots \times \mathcal X_n$ and 
can typically be done much more quickly than the first step. It can
be done in parallel since the projection onto $\mathcal X$ 
can be found by projections onto $\mathcal X_i$ for $i=1,\ldots,n$.
The third step is simply a dual update and is computationally 
inexpensive.
 
\subsection{Preconditioning}
Both theoretical analysis and practical evidence suggest that the 
precision and convergence rate of first-order methods can be
significantly improved by preconditioning the problem.
Here, we use 
diagonal scaling as preconditioning as discussed in 
\cite{beck2014introduction} and \cite{wright1999numerical}. 
Diagonal scaling can be viewed as applying an appropriate 
linear transformation before running the algorithm. When the set
$\mathcal X$ is convex, the preconditioning can substantially affect the
speed of convergence,
but does not affect the quality of the point returned,
(which must be a solution to the convex problem).
In other words, for convex problems,
preconditioning is simply a tool to help the algorithm converge faster. 
Optimal choice of preconditioners, even in the convex case, is still
an active research area 
\cite{giselsson2014diagonal,giselsson2014preconditioning,
giselsson2014improved, giselsson2014monotonicity,
ghadimi2015optimal, shi2014linear, hong2012linear,
boley2013local, deng2012global}. 
In the nonconvex case, however, preconditioning can have a critical role in 
the \emph{quality} of approximate solution, 
as well as the speed at which this solution is found.

Specifically, let $F\in\reals^{n\times n},E\in\reals^{m\times m}$
be diagonal matrices with positive diagonal entries. 
The goal is to choose $F$ and $E$ such that running ADMM on the 
following problem has better convergence properties 
\BEQ\label{prec}
\begin{array}{ll}
\mbox{minimize}   & (1/2)x^T P x + q^T x + I_{\mathcal X}(z)\\
\mbox{subject to} &  \left[ \begin{array}{c}
EA \\
F \end{array} \right]x - 
\left[ \begin{array}{c}
0 \\
F \end{array} \right]z = 
\left[ \begin{array}{c}
Eb \\
0 \end{array} \right]. \\
\end{array}
\EEQ

We use the choice of $E$ and $F$ recommended in 
\cite{giselsson2014diagonal} 
to minimize the effective condition number (the ratio of the largest 
singular value to the smallest non-zero singular value) of the following
matrix
\[
\left[ \begin{array}{cc}
E & 0\\
0 & F \end{array} \right]
\left[ \begin{array}{c}
A \\
I \end{array} \right]
P^\dagger 
\left[ \begin{array}{cc}
A^T &I \end{array} \right]
\left[ \begin{array}{cc}
E & 0\\
0 & F \end{array} \right],
\] 
where $P^\dagger$ denotes the pseudo-inverse of $P$.
Given matrix $M\in\reals^{n\times n}$, minimizing the condition 
number of $DMD$ for diagonal $D\in\reals^{n\times n}$ can be
cast as a semidefinite program. 
However, a heuristic called \emph{matrix equilibration} can be used to avoid the
computational cost of solving a semidefinite program.
(See \cite{sluis1969condition, bradley2010algorithms} and references 
therein.) Since for embedded applications computational 
resources are limited, we avoid finding $P^\dagger$ or equilibrating
completely. 
We instead find $E$ to normalize the rows of $A$ 
(usually in $\ell_1$ or $\ell_2$ norm) and set $F$ to be the identity.

After finding $E$ and $F$,
preconditioned ADMM has the following form:
\begin{equation}
\label{scaled_admm}
  \begin{split}
x^{k+1/2} &:= \left[ \begin{array}{cc}
I & 0\\
\end{array} \right]
\left[ \begin{array}{cc}
P+\rho F^2 & A^TE \\
EA & -(1/\rho)I 
\end{array} \right]^{-1}
\left[ \begin{array}{c}
-q+\rho \left(F^2x^k+A^TE^2b-
\left[ \begin{array}{cc}
A^TE & F 
\end{array} \right]
u^k\right)\\
0
\end{array} \right]\\
x^{k+1} &:=  \Pi\left(x^{k+1/2}+
\left[ \begin{array}{cc}
0 & F^{-1}
\end{array} \right]
u^k\right)\\
u^{k+1} &:=  u^k + \left[ \begin{array}{c}
EA \\
F \end{array} \right]x^{k+1/2} - 
\left[ \begin{array}{c}
0 \\
F \end{array} \right]x^k - 
\left[ \begin{array}{c}
Eb \\
0 \end{array} \right].
\end{split}
\end{equation}

\subsection {The overall algorithm}
We use the update rules \eqref{scaled_admm} for $k=1,\ldots,N$, where
$N$ denotes the (fixed) maximum number of iterations. Also, as 
described above, the algorithm is repeated for $M$ number of random
initializations. The computational cost of the algorithm consists of a 
factorization and $MN$ matrix products and projections. 
Here is a description of the overall algorithm with $f(x)=(1/2)x^T P x + q^T x + r$. 

\begin{minipage}[c]{1 \textwidth}
\begin{algorithm}[H]
\caption{Approximately solving nonconvex constraint QP~(\ref{miqp})}
\begin{algorithmic}
\State 
\If{$A$ or $P$ changed}
\State find $E$ and $F$ by equilibrating 
$\left[ \begin{array}{c}
A \\
I \end{array} \right]
P^\dagger 
\left[ \begin{array}{cc}
A^T &I \end{array} \right]$
\State find and store $LDL$ factorization of $\left[ \begin{array}{cc}
P+\rho F^2 & A^TE \\
EA & -(1/\rho)I 
\end{array} \right]$
\EndIf
\State $x_\mathrm{best} := \emptyset$, $f(x_\mathrm{best}) := \infty$ 
\For{random initialization $1, 2, \ldots, N$}
\For{iteration $1,2,\ldots,M$}
\State update $x$ from \eqref{scaled_admm}
\If {$\|Ax-b\|\leq\epsilon^{\mathrm{tol}}$ and $f(x)<f(x_{\textrm{best}})$}
\State $x_{\textrm{best}} = x$
\EndIf
\EndFor
\EndFor \\
\Return $x_{\textrm{best}}$.
\end{algorithmic}
\label{alg_summary}
\end{algorithm}
\end{minipage}
\bigskip

We mention a solution refinement technique here that can be used to 
find a solution with possibly better objective value after the algorithm stops. This 
technique, sometimes known as \emph{polishing} consists of fixing the 
nonconvex variable and solving the resulting convex optimization problem.
Using this technique, one may use larger $\epsilon^{\mathrm{tol}}$ during 
the $N$ iterations and only reduce $\epsilon^{\mathrm{tol}}$ at the refinement
step. Depending on the application, it might be computationally sensible to solve
the resulting convex optimization problem. 
Another effective technique is to introduce a notion of \emph{no-good cut} 
during iterations for problems with binary variables. 
A no-good cut prohibits the integer part to be equal to the previous one, by imposing 
one additional inequality constraint 
$\sum_{i\in T}x^{k+1/2}_{b_i} - \sum_{i\in F}x^{k+1/2}_{b_i} \leq B - 1$, where
$x_{b_1},x_{b_2},\ldots$ are binary variables and
$T=\{i| x^k_{b_i}=1\}, F=\{i| x^k_{b_i}=0\}$, and $B$ is the number of elements of 
$T$.
We do not use either of these techniques in 
the following examples.

\section{Numerical examples} 
\label{s-examples}
In this section, we explore the performance of our proposed algorithm
on some example problems. 
For each example, $\rho$ was chosen between $0.1$
and $10$ to yield good performance;
all other algorithm parameters were kept constant.
As a benchmark, we compare our results to the commercial 
solver MOSEK, which can globally solve MIQPs.
All experiments were carried out on
a system with two $3.06$ GHz cores with $4$ GB of RAM.

The results suggest that this heuristic is effective in finding 
approximate solutions for mixed integer quadratic programs. 
\subsection{Randomly generated QP}
First we demonstrate the performance of our algorithm
qualitatively for a random mixed-Boolean quadratic program.
The matrix $P$ in \eqref{miqp} was chosen 
as  $P=QQ^T$, where the entries of $Q\in \reals^{n\times n}$,
as well as those of $q$ and $A$,
were drawn from a standard normal distribution.
The constant $r$ was chosen such that the optimal value of the 
unconstrained quadratic minimization is $0$.
The vector $b$ was chosen
as $b=Ax_0$, where $x_0\in\mathcal X$ was chosen uniformly randomly,
thus ensuring that the problem is feasible.
We used $n=200$ and $m=50$ with $\mathcal X_i=\{0,1\}$ for $i=1,
\ldots,100$, $\mathcal X_i=\reals_+$ for $i=101,\ldots,150$, and
$\mathcal X_i=\reals$ for the other indices $i$.

We used MOSEK to find the optimal value for the problem. After 
$60832$ seconds (more than $16$ hours), MOSEK certifies that the optimal
value is equal to $2040$. 
We ran algorithm \ref{alg_summary} for $10$ different 
initializations and $200$ iterations for each initialization,
with step size $\rho=0.5$. For a naive implementation in MATLAB,
it took $120$ milliseconds 
to complete all precomputations (preconditioning and factorization),
and $800$ milliseconds to do all $2000$ iterations.
The best objective value found for the problem was $2067$ 
($1.3\%$ suboptimal).

One interesting observation is that the parameter $\rho$ tends to
trade off feasibility and optimality:
with small values of $\rho$, the algorithm often fails to find a feasible point,
but feasible points found tend to have low objective value. On the 
other hand, with large values of $\rho$, 
feasible points are found more quickly, but tend to have higher objective value.
\subsection{Hybrid vehicle control}
We consider a simple hybrid electric vehicle drivetrain
(similar to that of \cite[Exercise 4.65]{boyd2004convex}),
which consists of a battery, an electric motor/generator, and a heat engine,
in a parallel  configuration.
We assume that the demanded power $P^{\rm des}_t$
at the times $t = 0, \ldots, T-1$ is known in advance.
Our task is to plan out the battery and engine power outputs $P^{\rm batt}_t$
and $P^{\rm eng}_t$, for $t = 0, \ldots, T-1$,
so that
\[
P^{\rm batt}_t + P^{\rm eng}_t \geq P^{\rm des}_t.
\]
(Strict inequality above corresponds to braking.)

\paragraph{Battery.}
The battery has stored energy $E_t$ at time $t$,
which evolves according to 
\[
E_{t+1} = E_t - \tau P^{\rm batt}_t, \qquad t = 0, \ldots, T-1,
\]
where $\tau$ is the length of each discretized time interval.
The battery capacity is limited,
so that $0 \leq E_t \leq E^{\rm max}$ for all $t$,
and the initial energy $E_0$ is known.
We penalize the terminal energy state of the battery according to $g(E_T)$, 
where
\[
g(E) = \eta (E^{\rm max} - E)^2,
\]
for $\eta \geq 0$.

\paragraph{Engine.}
At time $t$, the engine may be on or off,
which is modeled with binary variable $z_t$.
If the engine is on ($z_t=1$),
then we have $0 \leq P^{\rm eng}_t \leq P^{\rm max}$,
and 
$\alpha (P^{\rm eng}_t)^2 + \beta P^{\rm eng}_t + \gamma$ units of fuel are consumed,
for nonnegative constants $\alpha$, $\beta$, and $\gamma$.
If the engine is off ($z_t = 0$), it consumes no fuel, and $P^{\rm eng}_t = 0$.
Because $z_t \in \{0,1\}$,
the power constraint can be written as $0 \leq P^{\rm eng} \leq P^{\rm max} z_t$,
and the fuel cost as $f(P^{\rm eng}_t, z_t)$, where
\[
f(P, z) = \alpha P^2 + \beta P + \gamma z.
\]
Additionally, we assume that turning the engine on after it has been off
incurs a cost $\delta \geq 0$,
\ie, at each time $t$, we pay $\delta (z_t - z_{t-1})_+$,
where $(\cdot)_+$ denotes the positive part.

\paragraph{Optimal power split problem.}
The hybrid vehicle control problem can be formulated as
\begin{equation}
\begin{array}{ll}
\mbox{minimize} & \eta (E_T - E^{\rm max})^2 + 
     \sum_{t=0}^{T-1} f(P^{\rm eng}_t, z_t) + \delta (z_{t} - z_{t-1})_+ \\
\mbox{subject to} 
  & E_{t+1} = E_t - \tau P^{\rm batt}_t \\
  & P^{\rm batt}_t + P^{\rm eng}_t \geq P^{\rm des}_t \\
  & z_t \in \{0, 1\},
\end{array}
\label{e-hybrid-opt-prob}
\end{equation}
where all constraints must hold for $t = 0, \ldots, T-1$.
The variables are 
$P^{\rm batt}_t$, $P^{\rm eng}_t$, and $z_t$ for $t = 0, \ldots, T-1$,
and $E_t$, for $t = 1, \ldots, T$.
In addition to the parameters given above,
we take $z_{-1}$ to be a parameter denoting the initial engine state.


\begin{figure}
\begin{center}
\begin{psfrags}
\psfrag{x}[B][B]{\raisebox{-1.2ex}{\tiny $t$}}
\psfrag{p1}[B][B]{\raisebox{+.5ex}{\tiny{$P^\mathrm{Eng}_t$}}}
\psfrag{p2}[B][B]{\raisebox{+.5ex}{{\tiny $P^\mathrm{batt}_t$}}}
\psfrag{p3}[B][B]{\raisebox{+.5ex}{\tiny{$E_t$}}}
\psfrag{p4}[B][B]{\raisebox{+.5ex}{\tiny{$z_t$}}}
\includegraphics[width=.46\linewidth]{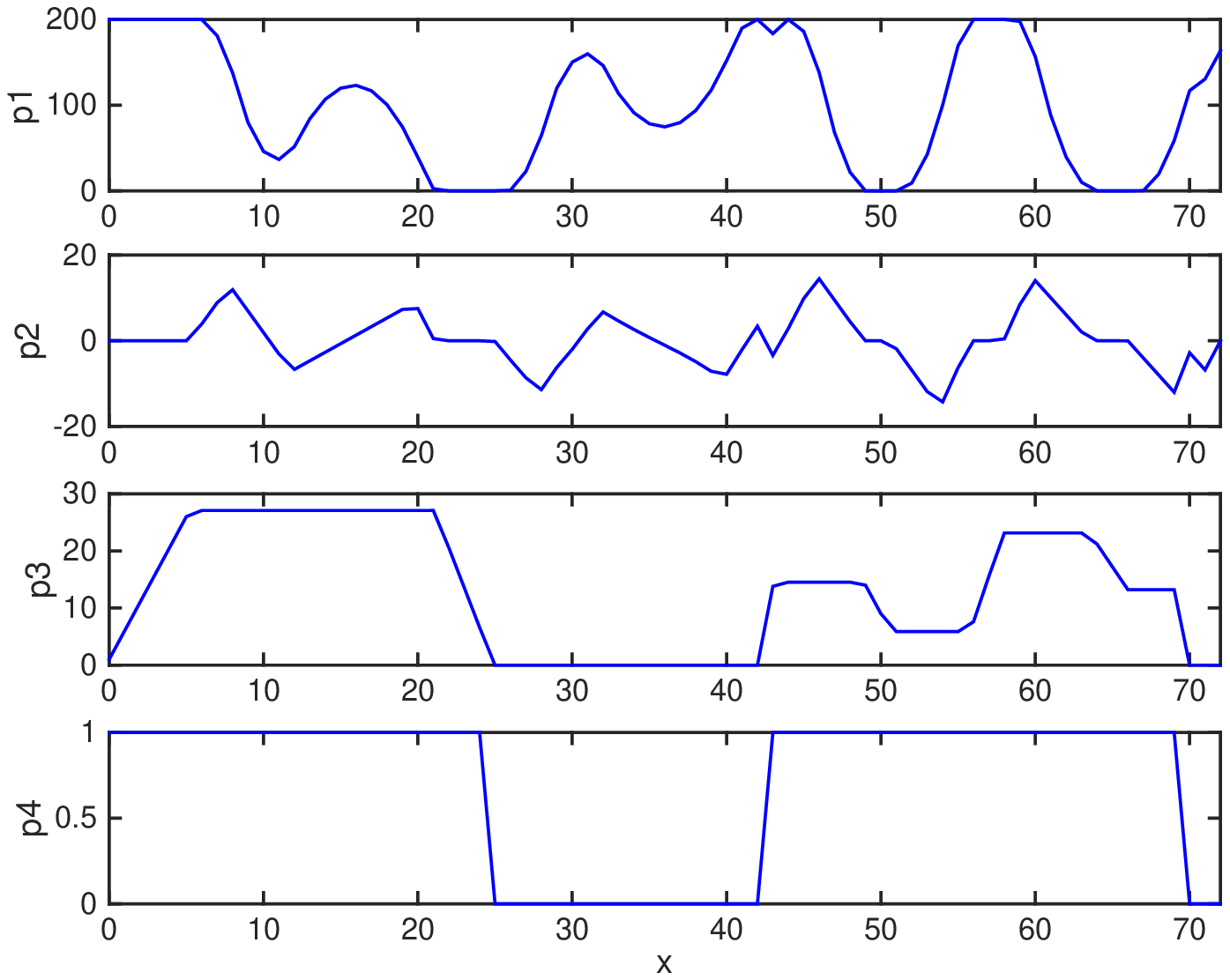} 
\end{psfrags}
\hspace*{\fill}
\begin{psfrags}
\psfrag{x}[B][B]{\raisebox{-1.2ex}{\tiny $t$}}
\psfrag{p1}[B][B]{\raisebox{+.5ex}{\tiny{$P^\mathrm{Eng}_t$}}}
\psfrag{p2}[B][B]{\raisebox{+.5ex}{{\tiny $P^\mathrm{batt}_t$}}}
\psfrag{p3}[B][B]{\raisebox{+.5ex}{\tiny{$E_t$}}}
\psfrag{p4}[B][B]{\raisebox{+.5ex}{\tiny{$z_t$}}}
\includegraphics[width=.46\linewidth]{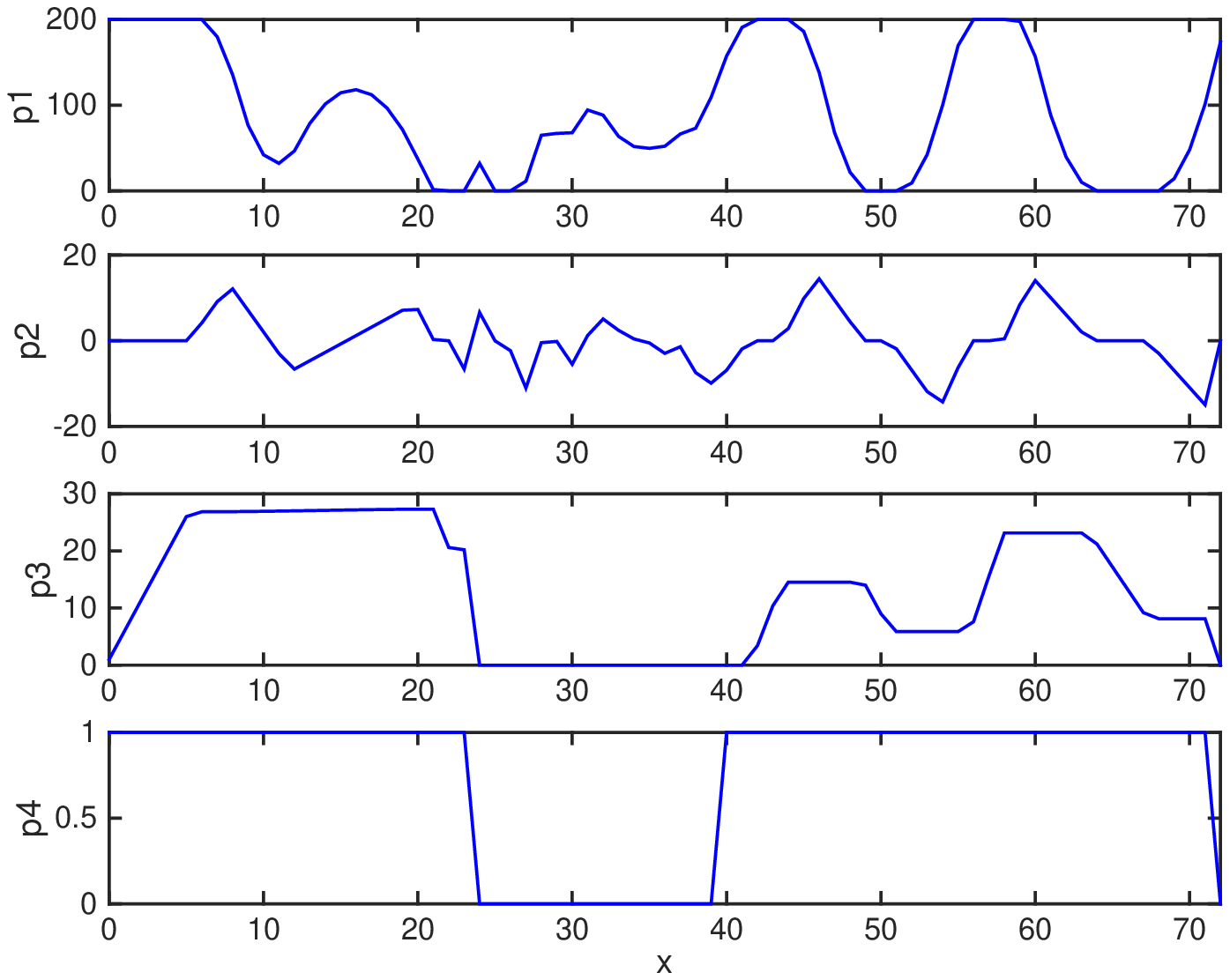} 
\end{psfrags}
\caption{
Engine power, battery power, battery energy, and engine on/off signals versus time.
Left: the global solution.
Right: the solution found using ADMM.
\label{vehicle}}
\end{center}
\end{figure}

We used the parameter values
$\alpha = 1$,
$\beta = 10$,
$\gamma = 1.5$,
$\delta = 10$,
$\eta = 0.1$,
$\tau = 5$,
$P^{\rm max} = 1$,
$E^{\rm max} = 200$,
$E_0 = 200$,
and 
$z_{-1} = 0$.
The demanded power trajectory $P^{\rm des}_t$
is not shown, but can be obtained by summing the engine
power and battery power in Figure \ref{vehicle}.
We ran the algorithm with $\rho=0.4$ for $1000$ iterations from
$5$ different initializations,
with primal optimality threshold $\epsilon^{\mathrm{tol}}=10^{-4}$.
The global solution found by MOSEK generates an objective
value of $339.2$ and the best objective value with our algorithm 
was $375.7$. 
In Figure \ref{vehicle},
we see that qualitatively, the optimal trajectory
and the trajectory generated by ADMM are very similar.

\subsection{Power converter control}

\begin{figure}
\centering
\resizebox{.5\textwidth}{!}{
\begin{tikzpicture}[scale=1.1]
\tikzstyle motor=[color=gray!50!white,draw=black,thick]
\ctikzset{bipoles/diode/height=.375}
\ctikzset{bipoles/diode/width=.3}

\draw[scale=1] (0,0) 
to [L, l=$L_1$] ++(3,0) node(N1){}
to ++(0,1)
to [L, l=$L_2$] ++(3,0) 
to ++(0,-1) node(N2){}
to ++(1.5,0) node(N3){}
to ++(2,0)
to [R, l=$R$] ++(0,-3)
-- (0, -3)
to [american voltage source, l=$u_t V_{\rm dc}$] ++(0,3)
;

\draw[scale=1] 
(N1)+(0,1) 
to ++(0,-1)
to [C, l=$C_1$] ++(3,0) 
to ++(0,1) node(N2){}
;

\draw[scale=1] 
(N3) to [C, l_=$C_2$] ++(0,-3)
;

\end{tikzpicture}
}
\caption{Converter circuit model.}
\label{f-converter}
\end{figure}
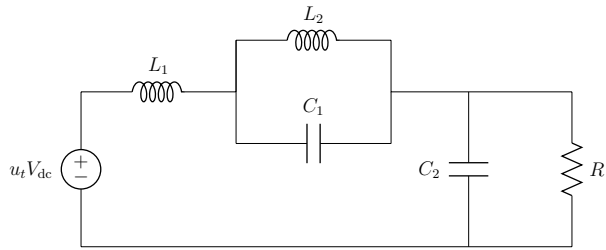

We consider control of the switched-mode power converter
shown in Figure \ref{f-converter}.
The circuit dynamics are
\[
\xi_{t+1} = G \xi_t + H u_t, \qquad t = 0, 1, \ldots, T-1,
\]
where $\xi_t = (i_{1,t}, v_{1,t}, i_{2,t}, v_{2,t})$ is the system state at epoch $t$,
consisting of all inductor currents and capacitor voltages,
and $u_t \in \{-1, 0, 1\}$ is the control input.
The dynamics matrices $G\in \reals^{4\times 4}$ and $H\in \reals^{4\times 1}$
are obtained by discretizing the dynamics of the circuit in Figure~\ref{f-converter}.

We would like to control the switch configurations
so that $v_2$ tracks a desired sinusoidal waveform.
This can be done by solving
\begin{equation}
\begin{array}{ll}
\mbox{minimize} & \sum_{t=0}^{T} (v_{2,t} - v_{\rm des})^2 
     + \lambda |u_{t} - u_{t-1}| \\
\mbox{subject to} 
  & \xi_{t+1} = G \xi_t + H u_t \\
  & \xi_0 = \xi_T \\
  & u_0 = u_T \\
  & u_t \in \{-1, 0, 1\},
\end{array}
\label{e-converter-opt-prob}
\end{equation}
where $\lambda \geq 0$ is a tradeoff parameter
between output voltage regulation and switching frequency.
The variables are $\xi_t$ for $t = 0, \ldots, T$
and $u_t$ for $t = 0, \ldots ,T-1$.

\begin{figure}
\begin{center}
\begin{psfrags}
\psfrag{t}[B][B]{\raisebox{-1.2ex}{\tiny $t$}}
\psfrag{u}[B][B]{\raisebox{0.5ex}{\tiny{$u_t$}}}
\psfrag{p}[B][B]{\raisebox{0.5ex}{\tiny{$v_{2,t}$}}}
\includegraphics[width=.46\linewidth]{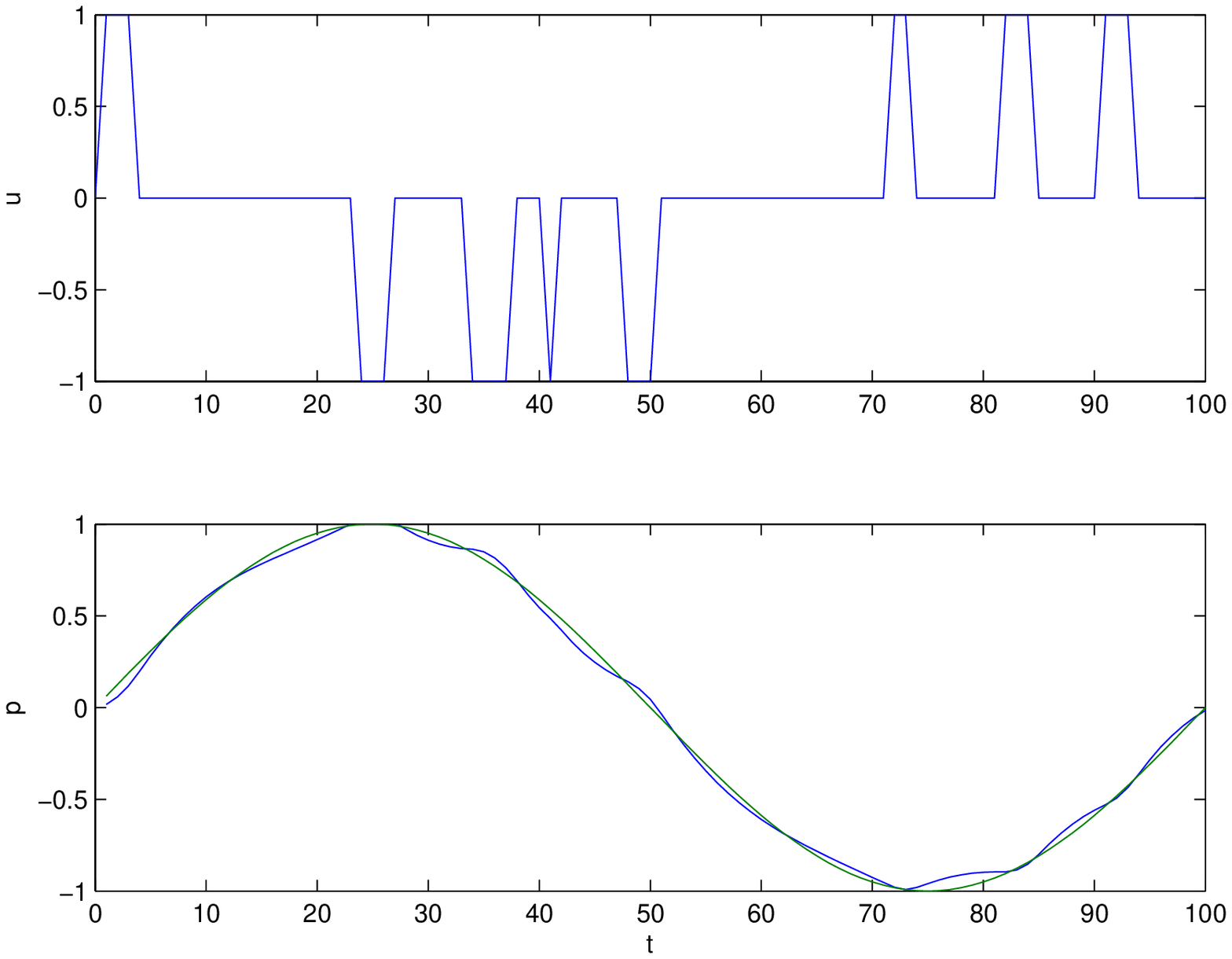} 
\end{psfrags}
\hspace*{\fill}
\begin{psfrags}
\psfrag{t}[B][B]{\raisebox{-1.2ex}{\tiny $t$}}
\psfrag{u}[B][B]{\raisebox{0.5ex}{\tiny{$u_t$}}}
\psfrag{p}[B][B]{\raisebox{0.5ex}{\tiny{$v_{2,t}$}}}
\includegraphics[width=.46\linewidth]{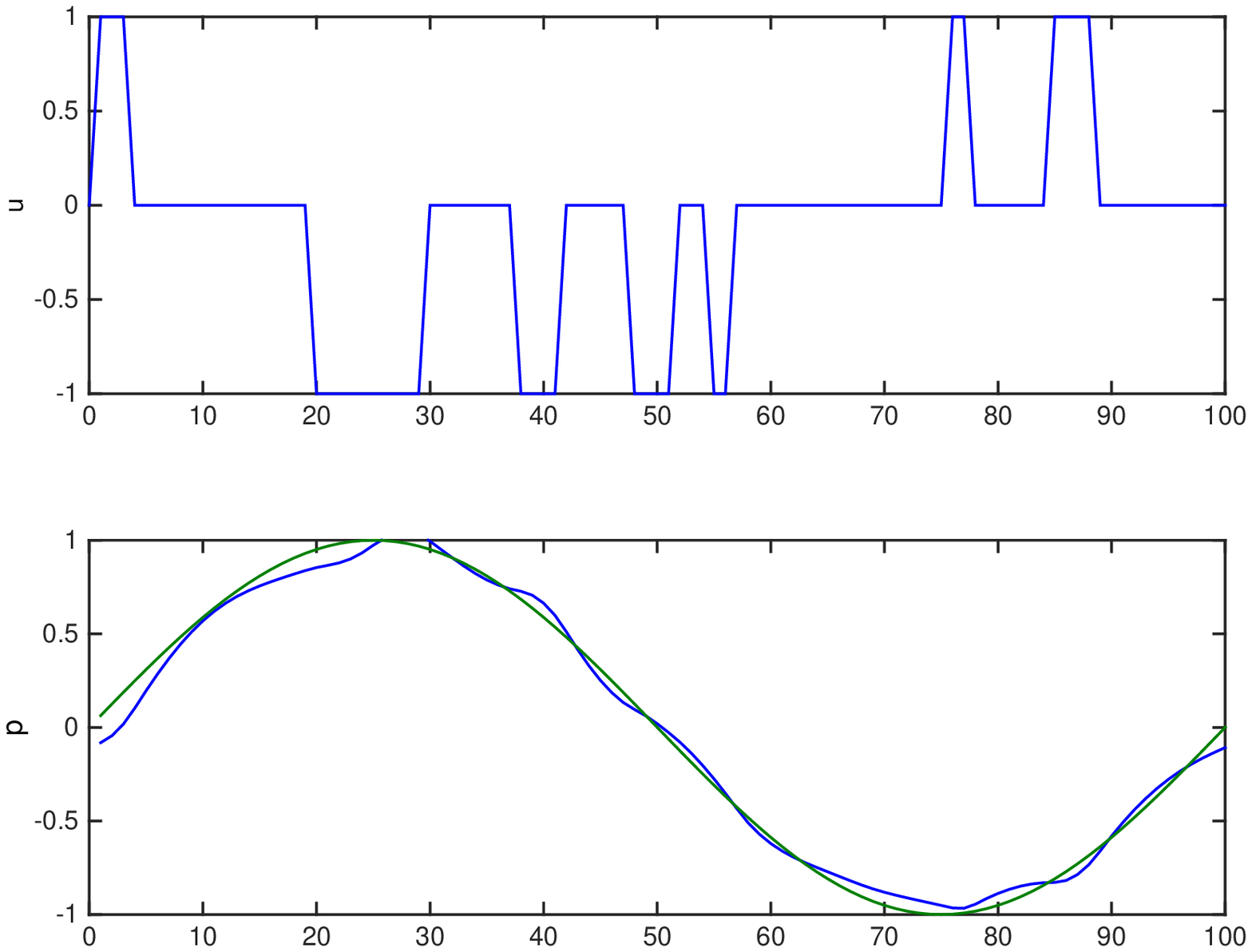} 
\end{psfrags}
\caption{The switch configuration and the output voltage. Left: the global solution. Right: the solution using ADMM.
\label{converter}}
\end{center}
\end{figure}

Note that if we take $\lambda = 0$, 
and take the input voltage $u_t$ to be unconstrained
(\ie, allow $u_t$ to take any values in $\reals$),
\eqref{e-converter-opt-prob} can be solved as a convex
quadratic minimization problem,
with solution $\xi_t^{\rm ls}$.
Returning to our original problem,
we can penalize deviation from this ideal waveform by including 
a regularization term $\mu \|\xi - \xi_t^{\rm ls}\|^2$ to \eqref{e-converter-opt-prob},
where $\mu>0$ is a positive weighting parameter.
We solved this regularized version of 
\eqref{e-converter-opt-prob}, with
$L_1=\SI{10}{\micro\henry}$,
$C_1 = \SI{1}{\micro \farad}$,
$L_2 = \SI{10}{\micro\henry}$,
$C_2 = \SI{10}{\micro \farad}$,
$R = \SI{1}{\ohm}$,
$V_{\rm dc} = \SI{10}{\volt}$,
$T = 100$ (with a discretization interval of $\SI{.5}{\micro\second}$),
$\lambda = \SI{1.5}{\volt^2}$,
and $\mu = 0.1$.
We run algorithm \ref{alg_summary} with $\rho=2.7$ and $500$ iterations for
three different initializations. An approximate solution is found via our heuristic
in less than $2$ seconds, whereas it takes MOSEK more than $4$ hours to find
the global solution.
Figure \ref{converter} compares the approximate solution derived by the heuristic
with the global solution.

\subsection{Signal decoding}
We consider maximum-likelihood decoding of a message 
passed through a linear multiple-input and multiple-output (MIMO) channel.
In particular, we have
\[
y = Hx + v,
\]
where
$y \in \reals^p$ is the message received,
$H \in \reals^{p\times n}$ is the channel matrix,
$x \in \reals^n$ is the message sent,
and the elements of the noise vector $v \in \reals^p$
are independent, identically distributed Gaussian random variables.
We further assume that the elements of $x$ belong to
the \emph{signal constellation} $\{-3, -1, 1, 3\}$.
The maximum likelihood estimate of $x$ is given by
the solution to the problem
\begin{equation}
\begin{array}{ll}
\mbox{minimize} & \| H\hat x - y \|^2 \\
\mbox{subject to} 
  & \hat x_i \in \mathcal \{-3, -1, 1, 3\}, \quad i = 1, \ldots, n,\\
\end{array}
\label{e-comms-opt-prob}
\end{equation}
where $\hat x \in \reals^n$ is the variable.

\begin{figure}
\begin{center}
\begin{psfrags}
\psfrag{x}[B][B]{\raisebox{-1.2ex}{\tiny $\|Hx^\textrm{rlx} - y \|_2^2 - \|Hx^\textrm{admm} - y\|_2^2$}}
\psfrag{y}[B][B]{\raisebox{0.5ex}{\tiny frequency}}
\includegraphics[width=.46\linewidth]{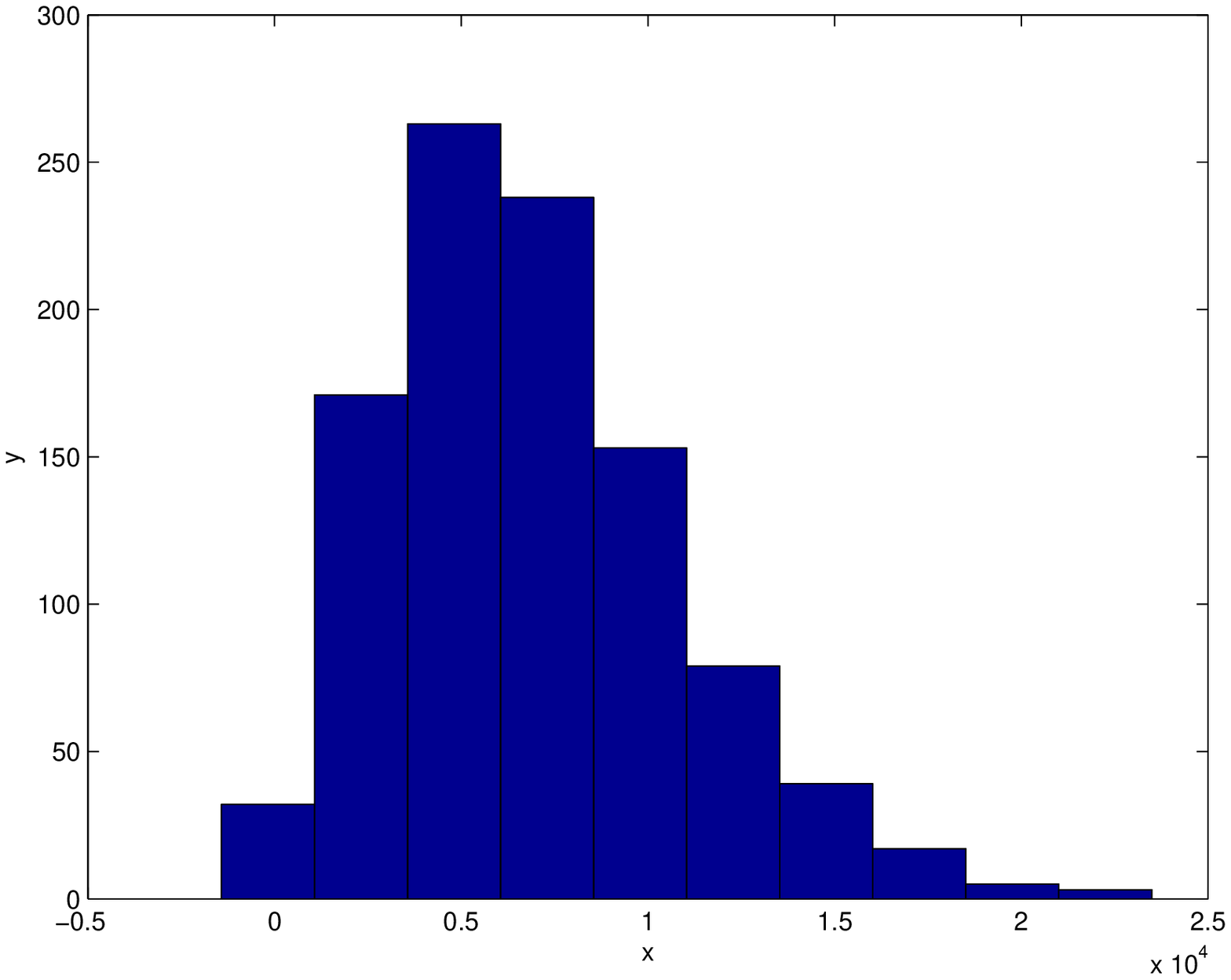} 
\end{psfrags}
\hspace*{\fill}
\begin{psfrags}
\psfrag{x}[B][B]{\raisebox{-1.2ex}{\tiny $\rm{BER}(x^\textrm{rlx}) - \rm{BER}(x^\textrm{admm})$}}
\psfrag{y}[B][B]{\raisebox{0.5ex}{\tiny frequency}}
\includegraphics[width=.46\linewidth]{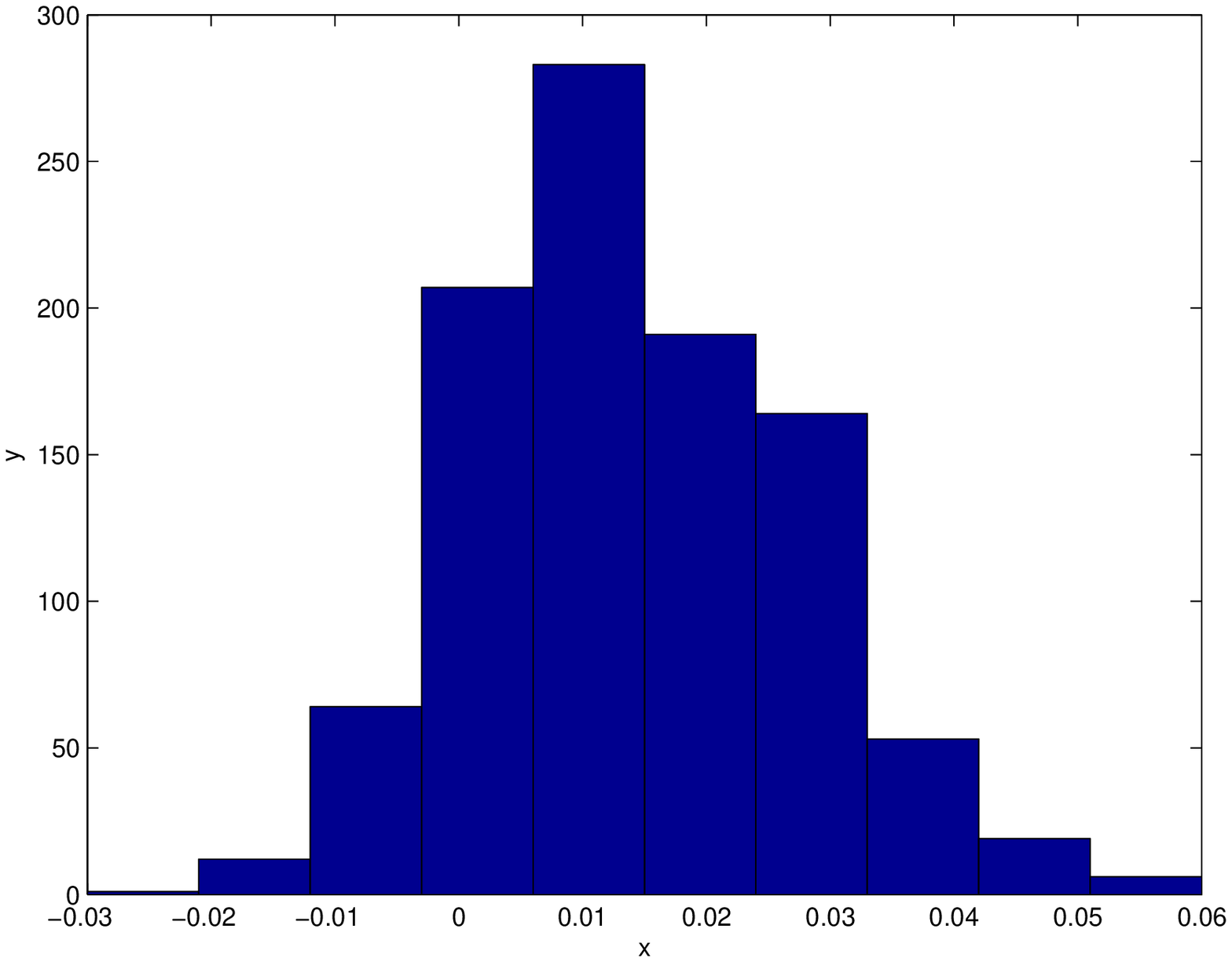} 
\end{psfrags}
\caption{Comparison of ADMM heuristic and relax-and-round.
Left: The difference in objective values. Right: The difference in bit error rates (BER).
\label{comms}}
\end{center}
\end{figure}

We generate $1000$ random problem instances with 
$H\in\reals^{2000\times400}$ 
chosen from a standard normal distribution. 
The uncorrupted signal $x$ is chosen 
uniformly randomly and the additive noise is Gaussian such that the 
signal to noise ratio (SNR) is $8$ dB.
For such a problem in embedded application, branch-and-bound methods 
are not desirable due to their worst-case time complexity.
We run the heuristic with only one initialization,
with $10$ iterations to find $x^\textrm{admm}$.
The average runtime for each problem (including preprocessing) 
is $80$ milliseconds, which is substantially faster than branch-and-bound 
based methods.
We compare the performance of the points $x^\textrm{admm}$
 with the points found by relax-and-round technique $x^\textrm{rlx}$.
 In Figure \ref{comms} we have plotted the histogram of the difference between the
 objective values evaluated at $x^\textrm{admm}$ and
 $x^\textrm{rlx}$. Depicted in Figure \ref{comms}, we see 
 that in $95\%$ of the cases,
the bit error rate (BER) using our heuristic was at least as good as the
bit error rate using relax and round. 

\section{Conclusions}
In this paper, we introduced an effective heuristic for finding approximate solutions 
to convex quadratic minimization problems over the intersection of affine and
nonconvex sets. Our heuristic is significantly faster than branch-and-bound
algorithms and has shown effective in a variety of embedded problems 
including hybrid vehicle control, power converter control, and signal decoding.

\nocite{*}
\bibliography{miqp}

\end{document}